\documentclass[12pt]{article}
\usepackage{amsmath, amssymb, amsthm}
\usepackage[margin=1.2in]{geometry}
\usepackage{times}

\newtheorem{defn0}{Definition}[section]
\newtheorem{prop0}[defn0]{Proposition}
\newtheorem{thm0}[defn0]{Theorem}
\newtheorem{lemma0}[defn0]{Lemma}
\newtheorem{corollary0}[defn0]{Corollary}
\newtheorem{example0}[defn0]{Example}
\newtheorem{remark0}[defn0]{Remark}
\newtheorem{conjecture0}[defn0]{Conjecture}

\newenvironment{definition}{ \begin{defn0}}{\end{defn0}}
\newenvironment{proposition}{\bigskip \begin{prop0}}{\end{prop0}}
\newenvironment{theorem}{\bigskip \begin{thm0}}{\end{thm0}}

\newcommand{\propref}[1]{Proposition~\ref{#1}}
\newcommand{\thmref}[1]{Theorem~\ref{#1}}

\newcommand{\secref}[1]{Section~\ref{#1}}

\newcommand{\m}{\mathfrak m}
\def\maxn{{\mathfrak n}}                   
\def\res{{\bf k}}                   

\def\soc{\operatorname{soc}}
\newcommand{\Hom}{\operatorname{Hom}}
\def\HF{\operatorname{HF}}

\def\rmod{{R\_mod}}
\def\rmodNoeth{{R\_{mod.Noeth}}}
\def\rmodArtin{{R\_mod.Artin}}

\title{  \bf \Huge {\sc{Inverse-syst.lib}} \\ Singular library for computing Macaulay's inverse systems
 }
\author{\large   Juan Elias
\thanks{Partially supported by  MTM2013-40775-P}
}

\date{\today }

\begin{document}

\maketitle

\begin{flushright}
\emph{
For J.L. S\'{a}nchez Palacios\\
who loved algebra with passion.
}
\end{flushright}

\medskip
\begin{abstract}
In this note we review  the Singular library {\sc{Inverse-syst.lib}}
that implements  Macaulay's correspondence and other related constructions for local rings.
\end{abstract}

\section{Introduction}

Let $\res$ be an arbitrary field.
Let $R=\res[\![x_1,\dots x_n]\!]$ be the ring of the formal series with maximal ideal $\m =(x_1,\cdots,x_n)$
and let  $S=\res[x_1,\dots,x_n]  $ be  a  polynomial ring, we  denote by $\m=(x_1,\dots,x_n)$ the homogeneous maximal ideal of $S$.

In 1916 Macaulay stated a one-to-one correspondence between Artin Gorenstein ideals of $R$  and polynomials of  $S$, \cite{Mac16}.
This correspondence can be extended to Artin ideals of $R$ and finitely sub-$R$-modules of $S$.
Recall that Macaulay's correspondence is a particular case of Matlis duality, see \thmref{matlis} and \propref{mac}.

Classically Macaulay's correspondence has been mainly used for studying homogenous ideals, \cite{Iar84}, \cite{Iar94}.
Recently Macaulay's correspondence has been applied to the classification of local Artin Gorenstein algebras,
see \cite{ER12},   \cite{CENR13}, \cite{ER15}, \cite{ES14}.
Most of the examples appearing in these papers have  been computed by using Singular, \cite{DGPS}.

In this note we review the main commands of the Singular library {\sc{Inverse-syst.lib}} that we used for these computations, \cite{E-InvSyst14}.
The main purpose of this library is to implement Macaulay's correspondence
if the action of $R$ in $E_R(\res)=S$ is defined by differentiation or contraction, \thmref{gab}.
We also implement some useful  operations of $S$ as $R$-module.
See \secref{comm} for a listing of all commands of {\sc{Inverse-syst.lib}}.
In \secref{el} we give a new proof of the classification of Artin Gorenstein local rings with Hilbert function
$\{1,3,3,1\}$ by using {\sc{Inverse-syst.lib}} and  the Weierstrass equation of an elliptic curve instead of Legendre equation  as it was done in \cite{ER12}.

\section{Macaulay's correspondence}

Let $A=R/I$ be an Artin ring with maximal ideal $\maxn=\m/I$.
The Hilbert function of $A$ is the numerical function
$\HF_A:\mathbb N \longrightarrow \mathbb N$ defined by
$\HF_A(i)=\dim_{\res}(\maxn^i/\maxn^{i+1})$, $i\ge 0$.
The socle degree of $A$ is the last integer $s$ such that $\HF_A(s)\neq 0$.
The socle of $A$ is the $\res$-vector subspace of $A$
$\soc(A)=(0:_A\maxn)$, and the  Cohen-Macaulay type of  $A$
is $\tau(A)=\dim_{\res}(\soc(A))$.

\begin{definition}\label{deflevel}An Artin ring  $A$ with socle degree $s$ is  level  if $\soc(A)=\maxn^s$.
In particular, $A$ is Gorenstein iff $\tau(A)=1$.
\end{definition}

The function {\tt{isAG(I)}}  returns $-2$ if the quotient $A=R/I$ is not Artin, returns $-1$ if $A$ is Artin but not
Gorenstein and returns the socle degree of $A$ if the ring $A$ is an Artin  Gorenstein ring.
A similar function is implemented  for checking if $A$ is level, see \secref{comm}.
The functions {\tt socle} and {\tt cmType} compute the socle and the Cohen-Macaulay type of $A$.

{\footnotesize{
\begin{verbatim}
// we define the ring r
>ring r=0,(x(1..3)),ds;
// loading the libary
>LIB "inverse-syst-v.4.lib";
// the quotient r/i is not Artin:
>ideal i=x(1)^2+x(2)^3, x(2)^4;
>isAG(i);
-2
// the quotient r/i is Artin but Gorenstein:
>ideal i=x(1)^2+x(2)^3, x(2)^4+x(1)^2, x(3)^2+x(1)*x(2),
    x(1)*x(2)^2*x(3);
>isAG(i);
-1
// the  quotient r/i is Artin Gorenstein and we get
// the socle degree
> ideal i=x(1)^2+x(2)^3, x(2)^4+x(1)^2, x(3)^2+x(1)*x(2);
> isAG(i);
4
// we define an Artin no Gorenstein ideal:
> ideal i=x(1)^2+x(2)^3, x(2)^4+x(1)^2, x(3)^2+x(1)*x(2),
    x(1)*x(2)^2*x(3);
> isAG(i);
-1
// we compute the socle ideal of r/i
> socle(i);
_[1]=x(1)^2
_[2]=x(1)*x(2)+x(3)^2
_[3]=x(2)^3
_[4]=x(2)^2*x(3)
_[5]=x(1)*x(3)^2
_[6]=x(2)*x(3)^2
_[7]=x(3)^3
// we compute the Cohen-Macaulay type of r/i
> cmType(i);
3
\end{verbatim}
}}

The polynomial ring $S$ can be considered as $R$-module with two linear structures: by derivation and by contraction.
If  $char(\res)=0$, the $R$-module structure of $S$ by derivation is defined by
$$
\begin{array}{cccc}
  R\times S & \longrightarrow & S &   \\
  (x^\alpha, x^\beta) & \mapsto & x^\alpha\circ x^\beta = &
  \left\{
    \begin{array}{ll}
      \frac{\beta!}{(\beta-\alpha)!}x^{\beta-\alpha} &  \beta \ge \alpha \\
      0, & {\text{otherwise}}
    \end{array}
  \right.
\end{array}
$$
where for all $\alpha, \beta\in\mathbb N^n$, $\alpha!=\prod_{i=1}^n \alpha_i!$

\noindent
If $char(\res)\ge 0$, the $R$-module structure of $S$ by contraction is defined by:
$$
\begin{array}{cccc}
  R\times S & \longrightarrow & S &   \\
  (x^\alpha, x^\beta) & \mapsto & x^\alpha \circ x^\beta = &
  \left\{
    \begin{array}{ll}
      x^{\beta-\alpha} &  \beta \ge \alpha \\
      0, & otherwise
    \end{array}
  \right.
\end{array}
$$

\bigskip
\noindent
In Singular we can use the above products  as follows:

{\footnotesize{
\begin{verbatim}
> ideal F=x(1)^2*x(3)^4+ x(2)^3*x(1)*x(3)+x(2)^5;
> diff(x(1)^2,F);
_[1,1]=2*x(3)^4
> contract(x(1)^2,F);
_[1,1]=x(3)^4
\end{verbatim}
}}

\noindent
It is easy to prove that for any field $\res$ there is a $R$-module homomorphism
$$
\begin{array}{cccc}
 \sigma : &(S,der) & \longrightarrow & (S,cont) \\
  &x^\alpha& \mapsto& \alpha! \; x^\alpha
\end{array}
$$
If $char(\res)=0$ then $\sigma$ is an isomorphism of $R$-modules.
The $R$-module $S$ is the injective hull $E_R(\res)$ of the $R$-module $\res$:

\begin{theorem}
(\cite{Gab59}) \label{gab}
If $\res$ is of characteristic zero then
$$
E_R(\res)\cong (S, der) \cong (S,cont).
$$
If $\res$ is of positive  characteristic  then
$$
E_R(\res)\cong  (S,cont).
$$
\end{theorem}

\medskip

Given a commutative ring $R$ we denote by $\rmod$, resp. $\rmodNoeth$, $\rmodArtin$, the category of $R$-modules, resp. category of Noetherian $R$-modules, Artinian $R$-modules.
The Matlis dual of an $R$-module $M$ is $M^{\vee}=\Hom_R(M,E_R(\res))$.
We  write $(-)^{\vee}=\Hom_R(-,E_R(\res))$, which is an additive  contravariant exact functor from the category $\rmod$ into itself.

\begin{theorem}[Matlis duality]
\label{matlis}
The functor $(-)^\vee$ is  contravariant, additive and exact, and defines  anti-equivalence between $\rmodNoeth$ and  $\rmodArtin$ (resp. between   $\rmodArtin$ and $\rmodNoeth$).
 The composition $(-)^\vee \circ (-)^\vee$  is the identity functor of $\rmodNoeth$ (resp. $\rmodArtin$).
Furthermore, if $M$ is a $R$-module of finite length then $\ell_R(M^\vee)=\ell_R(M)$.
\end{theorem}

From the previous result  we can recover the classical result of Macaulay, \cite{Mac16}, for the power series ring,
see \cite{Ems78}, \cite{Iar94}.
If $I\subset R$ is an  ideal, then $(R/I)^\vee$ is the  sub-R-module of $S$
$$
 {I^{\perp}}=\{ g \in S\ |\  I \circ g = 0    \},
$$
\noindent
this is the Macaulay's inverse system of $I$.
Given a sub-$R$-module $M$ of $S$ then dual $M^\vee$ is an ideal of $R$
$$ M^\perp= \{ f \in R \ \mid \ f \circ g= 0 \ \mbox{ for all \ } g \in M\}.
$$

\begin{proposition}[Macaulay's duality]
\label{mac}
Let  $R=\res[\![x_1,\dots x_n]\!]$ be the $n$-dimensional power series ring over a field $\res$.
There is a order-reversing bijection $\perp$ between the set of finitely generated sub-$R$-submodules of $S=\res[\![y_1,\dots y_n]\!]$ and the set of $\m$-primary ideals of $R$ given by:
if $M$ is a submodule of $S$ then $M^\perp=(0:_R M)$, and $I^\perp=(0:_S I)$ for an ideal $I\subset R$.
 \end{proposition}

Given a polynomial $H\in S$ of degree $l$ we denote by $\mbox{top}(H)$ the degree $l$ homogeneous form of $H$.

\begin{proposition}[Proposition 3.7 and Corollary 3.8 \cite{Iar84}, \cite{DeS12}]
\label{levelchar}
Let $A=R/I$ be an Artin ring of socle degree $s$ and Cohen-Macaulay $t$.
The following conditions are equivalent:
\begin{enumerate}
\item[(i)] $A$ is level,
\item[(ii)] $I^{\perp}$ is generated by $t$ polynomials $H_1,\ldots,H_t\in S$ such that $\deg(H_i)=s$, for $i=1,\ldots,t$, and the homogeneous forms  $\rm{top}(H_1),\ldots,\rm{top}(H_t)$ are $\res-$linear independent.
\end{enumerate}
In particular,  $A=R/I$ is Gorenstein of socle degree $s$ if and only if
$I^\perp$ is a cyclic $R$-module generated by a polynomial of degree $s$.
\end{proposition}

Given a collection of polynomials $H_1,\ldots,H_t\in S$ we denote by
$\langle H_1,\ldots,H_t \rangle$ the sub-$R$-module of $S$ generated by
$H_1,\ldots,H_t$.
Notice that $\langle H_1,\ldots,H_t \rangle$ is not an ideal of $S$, is the $\res$-vector space generated
by the collection $H_1,\ldots,H_t$ and their derivatives of any order.
Notice that in Singular ideals are handled by the list of a given system of generators.
In the library {\sc{Inverse-syst.lib}} the sub-$R$-modules of $S$ are handled  by using the
Singular's structure of ideal, i.e. by the list of a given system of generators.

We denote by $S_{\le i}$ (resp. $S_{< i}$, resp. $S_i$), $i\in \mathbb N$,  the $\res$-vector space of polynomials of $S$ of degree less or equal (resp. less, resp. equal to) to $i$, and we consider the following $\res$-vector space
$$
 (I^{\perp})_i := {\frac{ I^{\perp} \cap S_{\le i} +  S_{< i}}{ S_{< i}}}.
$$

\begin{proposition} (\cite{Eviasm})
\label{hf}
For all $i\ge 0$ it holds
$$
\HF_{A}(i) = \dim_{\res}   (I^{\perp})_i.
$$
\end{proposition}

\bigskip

In the library {\sc{Inverse-syst.lib}} Macaulay's correspondence has been programmed with respect the two
$R$-module structures of $S$. i.e. with respect to the differentiation and with respect to the contraction.
Here we will show how works Macaulay's correspondence with respect the differentiation.
Recall  that for technical reasons the sub-$R$-modules of $S$ are handled in {\sc{Inverse-syst.lib}}  by using the Singular's structure of ideal.

The command {\tt{invSyst}} computes the inverse system $I^\perp \subset S$ of an ideal $I$ of $R$, the command {\tt{idealAnn}} computes the annihilator
 $M^\perp\subset R$ of  a finitely sub-$R$-module $M$ of S.

In the next example we will show that the composition {\tt{idealAnn}} $\circ$ {\tt{invSyst}} is the identity map  on the set of Artin ideals as \propref{mac} predicts.

{\footnotesize{
\begin{verbatim}
> ring r=0,(x(1..3)),ds;
> LIB "inverse-syst-v.4.lib";
// we define an ideal i of r, notice that the first generator is
// a random polynomial with monomials of degree between 2 and 3
// and random coefficients between -2 and 2. The second generator
// is a random homogenous form of degree 3 and coefficients
// between -1 and 1.
> ideal i=genPol(2,3,2), genPol(3,3,1), x(2)^3+x(1)*x(3)^4,
    x(1)^2+x(2)^2*x(3);
> i;
i[1]=2*x(1)^2+2*x(2)^2-x(1)*x(3)+2*x(2)*x(3)-x(3)^2-2*x(1)^3
    +x(1)^2*x(2)+2*x(1)*x(2)^2-2*x(2)^3-2*x(1)^2*x(3)
    +2*x(1)*x(2)*x(3)+2*x(2)^2*x(3)-2*x(2)*x(3)^2-x(3)^3
i[2]=-x(1)^2*x(2)-x(2)^3+x(1)*x(2)*x(3)+x(2)^2*x(3)+x(1)*x(3)^2
    +x(3)^3
i[3]=x(2)^3+x(1)*x(3)^4
i[4]=x(1)^2+x(2)^2*x(3)
// we compute the inverse system of i:
> ideal iv=invSyst(i);
> iv;
iv[1]=3*x(1)^2+69*x(2)^2-42*x(1)*x(2)*x(3)-3*x(2)^2*x(3)
    -42*x(1)*x(3)^2+15*x(2)*x(3)^2+22*x(3)^3
iv[2]=24*x(1)*x(3)+3*x(1)*x(2)^2+6*x(1)*x(3)^2-2*x(3)^3
// Notice that iv is not cyclic so i is not Gorenstein
// we compute the annihilator of iv. Should be the ideal i.
> ideal j=idealAnn(iv);
> j;
j[1]=7*x(1)^2+x(1)*x(2)*x(3)
j[2]=14*x(2)^2-7*x(1)*x(3)+14*x(2)*x(3)-7*x(3)^2+
    28*x(1)*x(2)^2-52*x(1)*x(2)*x(3)+196*x(2)^2*x(3)
    -98*x(2)*x(3)^2
j[3]=343*x(1)*x(2)*x(3)-686*x(2)^2*x(3)+343*x(2)*x(3)^2
j[4]=5*x(1)*x(3)^2-8*x(2)*x(3)^2+5*x(3)^3
j[5]=x(2)*x(3)^3
j[6]=x(3)^4
// we check that i and j are the same ideal:
> eqIdeal(i,j);
1
\end{verbatim}
}}

In the next example we will show that the composition  {\tt{invSyst}} $\circ$ {\tt{idealAnn}}  is the identity
map on the set of finitely sub-$R$-modules of $S$ as \propref{mac} predicts.

{\footnotesize{
\begin{verbatim}
// We start with a random polynomial
> ideal q=genPol(2,3,2);
> q;
q[1]=2*x(1)^2-2*x(1)*x(2)+2*x(2)^2+2*x(1)*x(3)-2*x(2)*x(3)
    -x(3)^2-x(1)^3-2*x(1)^2*x(2)+2*x(1)*x(2)^2-2*x(2)^3
    -2*x(1)*x(2)*x(3)+x(2)^2*x(3)-x(2)*x(3)^2
// we compute the annihilator qa of q.
// The ideal qa is a Gorenstein ideal.
> ideal qa=idealAnn(q);
> qa;
qa[1]=4*x(1)^2-17*x(1)*x(2)-5*x(2)^2+12*x(2)*x(3)+x(1)^3
qa[2]=2*x(1)*x(2)+2*x(2)^2-8*x(1)*x(3)-3*x(1)^2*x(2)
qa[3]=x(1)*x(3)-x(3)^2+2*x(1)*x(2)*x(3)
qa[4]=x(2)^3-6*x(1)*x(2)*x(3)
qa[5]=x(2)^2*x(3)+x(2)*x(3)^2
qa[6]=x(3)^3
//  We get that qa is a Gorenstein ideal and the socle degree
//  of r/i is three that coincides with the degree of q:
> isAG(qa);
3
// we compute the inverse system of qa
> ideal q2=invSyst(qa);
> q2;
q2[1]=6*x(1)*x(2)-24*x(1)*x(3)+22*x(3)^2+17*x(1)^3
    +34*x(1)^2*x(2)-34*x(1)*x(2)^2+34*x(2)^3
    +34*x(1)*x(2)*x(3)-17*x(2)^2*x(3)+17*x(2)*x(3)^2
// from Macaulay's correspondence q and q2 should coincide:
> eqModIH(q,q2);
1
\end{verbatim}
}}

\section{A case study: Artin Gorenstein rings with Hilbert \\function $\{1,3,3,1\}$}
\label{el}

As a corollary of the main result of \cite{ER12} we got the classification of Artin Gorenstein local rings with
Hilbert function $\{1,3,3,1\}$ by using the Legendre equation of an elliptic curve.
In this section we recover this classification by using the Weierstrass equation of an elliptic curve and the library
{\sc{Inverse-syst.lib}}.

\begin{theorem} [\cite{ER12}] \label{V(F)}
The classification of   Artin  Gorenstein local $\res$-algebras
 with Hilbert function  $\HF_A=\{1,n,n,1\}$ is equivalent
to the projective classification of the  hypersurfaces
$ V(F)\subset \mathbb {P}^{n-1}_{\res}$ where  $F$ is  a degree three
non degenerate form in $n$ variables.
\end{theorem}

See
\cite{ER12} for the classification of Artin  Gorenstein local $\res$-algebras
 with Hilbert function  $\HF_A=\{1,n,n,1\}$ for $n=1,2$.
Assume $n=3$.
 Any plane elliptic cubic curve $C\subset \mathbb P^{2}_{\res}$ is defined, in a suitable
system of coordinates, by a   Weierstrass' equation, \cite{Sil09} proof of Proposition 1.4,
$$
W_{a,b})=x_2^2 x_3=x_1^3+a x_1 x_3^2+ b x_3^3
$$
with $a,b \in \res$ such that  $4a^3+27b^2\neq 0$.
The $j$ invariant of $C$ is
$$
j(a,b)=1728\frac{4a^3}{4a^3+27b^2}
$$
It is well known that
two plane elliptic  cubic  curves $C_i=V(W_{a_i,b_i})\subset \mathbb P^{2}_{\res}$, $i=1,2$, are projectively isomorphic
if and only if $j(a_1,b_1)=j(a_2,b_2)$.

For elliptic curves the inverse moduli problem can be  done as follows.
We denote by $W(j)$ the following elliptic curves with $j$ as moduli :
$W(0)=x_2^2x_3+x_2x_3^2-x_1^3$,
$ W(1728)=x_2^2x_3-x_1x_3^2-x_1^3$,
and  for $j\neq 0, 1728$
$$
W(j)=(j-1728) (x_2^2 x_3+x_1 x_2 x_3 -x_1^3) + 36 x_1 x_3^2 +x_3^3.
$$

We will show by using the library
{\sc{Inverse-syst.lib}} that:

\begin{proposition}
Let $A$ be an Artin   Gorenstein local $\res$-algebra   with Hilbert function
$\HF_A=\{1,3,3,1\}$. Then $A$ is isomorphic to one and only one of the following
  quotients of  $R=\res[[x_1,x_2,x_3]]$ :
$$
\hskip -.1cm
\begin{footnotesize}
\begin{array}{|c|c|c|}           \hline
 \text{ Model for } A=R/I& \text{Inverse system } F & \text{Geometry of } C=V(F)\subset \mathbb P^{2}_{\res} \\ \hline
 (x_1^2,x_2^2,x_3^2) & x_1 x_2 x_3 & \text{Three independent lines} \\  \hline
 (x_1^2,x_1x_3,x_3x_2^2,x_2^3,x_3^2+x_1x_2) & x_2(x_1x_2- x_3^2) & \text{Conic and a tangent line} \\  \hline
 (x_1^2,x_2^2,x_3^2+6x_1x_2) &  x_3(x_1x_2- x_3^2)& \text{Conic and a non-tangent line} \\  \hline
 (x_3^2, x_1x_2,x_1^2+x_2^2-3 x_1x_3)&  x_2^2 x_3- x_1^2( x_1+ x_3)& \text{Irreducible nodal cubic} \\  \hline
 (x_3^2,x_1x_2, x_1x_3,x_2^3,x_1^3+ 3x_2^2x_3) & x_2^2 x_3 - x_1^3 & \text{Irreducible cuspidal cubic} \\  \hline
  (x_3^3, x_1^3+3x_2^2x_3, x_1x_3, x_2^2-x_2x_3+x_3^2, x_1x_2) & W(0)   &  \text{Elliptic curve $j=0$} \\ \hline
(x_2^2+x_1x_3, x_1x_2, x_1^2-3x_3^2) & W(1728)  &  \text{Elliptic curve $j=1728$} \\ \hline
 I(j)=(x_2(x_2- 2x_1),H_j,G_j) & W(j)& \text{Elliptic curve with $j \neq 0,1728$} \\ \hline
 \end{array}
 \end{footnotesize}
 $$
 with:

 $H_j=6jx_1x_2- 144(j-1728)x_1x_3+72(j-1728)x_2x_3- (j-1728)^2 x_3^2$, and

$G_j=jx_1^2-12(j-1728)x_1x_3+ 6(j-1728)x_2x_3+ 144(j-1728)x_3^2$;

 \noindent $I(j_1)\cong  I(j_2)$ if and only if $j_1=j_2$.
\end{proposition}

\medskip
The first $7$ models can be obtained from the corresponding inverse system $F$
by using the command {\tt{idealAnn}}.
Assume that $j\neq 0, 1728$.
Let $J(j)$ be the ideal $\langle W(j)\rangle^\perp$;
a simple computation shows that $\HF_{R/J(j)}=\{1,3,3,1\}$, \propref{hf}.

{\footnotesize{
\begin{verbatim}
// we define a ring  of characteristic zero, three variables and ground field
// a field of functions with indeterminate c(1)
> def r=workringc(0,1,3);
> setring r;
> r;
//   characteristic : 0
//   1 parameter    : c(1)
//   minpoly        : 0
//   number of vars : 3
//        block   1 : ordering ds
//                  : names    x(1) x(2) x(3)
//        block   2 : ordering C
// the ideal p defines a elliptic curve with j=c(1)
> ideal p=weierstrassp();
> p;
p[1]=(-c(1)+1728)*x(1)^3+(c(1)-1728)*x(1)*x(2)*x(3)+
    (c(1)-1728)*x(2)^2*x(3)+36*x(1)*x(3)^2+x(3)^3
// we define the ideal q. We will prove that the inverse system of q
// is p.
> ideal q=idealwp();
> q;
q[1]=(6*c(1))*x(1)*x(2)+(-144*c(1)+248832)*x(1)*x(3)+
    (72*c(1)-124416)*x(2)*x(3)+(-c(1)^2+3456*c(1)-2985984)*x(3)^2
q[2]=(c(1))*x(1)^2+(-12*c(1)+20736)*x(1)*x(3)+(6*c(1)-10368)*x(2)*x(3)+
    (144*c(1)-248832)*x(3)^2
q[3]=-2*x(1)*x(2)+x(2)^2
// we check that q is contained in p^\perp
> diff(q,p);
_[1,1]=0
_[1,2]=0
_[1,3]=0
// If we perform de division of the $4$-th power of the maximal ideal
// by q we get three matrices Q, R, U such that (see Singular's manual)
//           generators(maxideal(4))*U=generators(q)*Q + R
// U is the 15x15 identity matrix, R is the 15x1 zero matrix
// and Q is a 6x15 matrix with coefficients in the ground field (see below
// for more details). The command is:
> division(maxideal(4),q);
\end{verbatim}
}}

\noindent
From the last computation we get that
the denominators of the coefficients of the matrix $Q$ are constant polynomials or polynomials
with roots in $\{0, 1728\}$.
Hence for all $j=c(1)\neq 0, 1728$ we get that $m^4\subset q$, so $q$ is an Artin ideal.
Notice that for all $j=c(1)\in \res$,  $q=I(j)$ and $p=\langle W(j)\rangle$.
Since $I(j)$ is generated by three homogeneous elements, $I(j)$ is a homogeneous complete intersection ideal.
In particular $I(j)$ is a homogeneous Artin Gorenstein ideal, so $\HF_{R/I(j)}$ is symmetric.
Notice that the generators of $I(j)$ are three homogeneous forms of degree two, so the Hilbert function
of $A=R/I(j)$ is $\{1,3,3,1\}$.
We know that $I(j)$ is contained $J(j)=\langle W(j)\rangle^\perp$.
Since $\HF_{R/I(j)}=\HF_{R/J(j)}=\{1,3,3,1\}$, we get that $I(j)=J(j)=\langle W(j)\rangle^\perp$,
i.e. $I(j)=\langle W(j)\rangle^\perp$.

\section{Commands}
\label{comm}

Next,  we list the most important commands of {\sc{Inverse-syst.lib}}.

{\footnotesize{
\begin{verbatim}

IDEAL THEORY
genPol(i,j,a).
    Returns a generic polinomial sum of forms of degrees 	
    between i and j, with integer coefficients in [-a,a].
eqIdeal(J,I); I, J ideals.
    Returns 1 if I=J, 0 otherwise.
socle(J); I ideal
    Returns -1 if J is not Artin,
    returns the  ideal of J if J is Artin.
cmType(J); J ideal.
    Returns -1 if J is not Artin,
    returns the  Cohen-Macaulay type of J otherwise.
isAG(I); I ideal.
    Returns -2 if J is not Artin,
    returns -1 if J is 	Artin but not Gorenstein, 	and
    returns the socle degree if J 	is Artin Gorenstein.
isLevel(I); I ideal.
    Returns -2 if J is not Artin,
    returns -1 if is Artin but not Level, and
    returns the socle degree if J is Artin Level

MACAULAY INVERSE SYSTEM CORRESPONDENCE WITH COEFFICIENTS

invSystG(ideal J)
    returns the inverse system of J; J is Artin Gorenstein
invSyst(J)
    returns the inverse system of J; J is Artin
idealAnnG(poly f)
    returns the Artin Gorenstein ideal with inverse system f
idealAnn(I)
    returns the Artin  ideal with inverse system I

STRUCTURE OF INJECTIVE HULL WITH COEFFICIENTS

memberIH(g,I); I=f1,... list of polynomials, g polynomial.
    returns 1 if g belongs to the R-submodule of S generated by f1,...fs in S,
    returns 0, otherwise
subModIH(I,J); I=f1,... list of polynomials, J=g1,... list of polynomials.
    Returns 1 if I is a sub-R-submodule of J, both sub-R-modules of S;
    0 otherwise
eqModIH(I,J); I=f1,... list of polynomials, J=g1,... list of polynomials.
    Returns 1 if I=J, both sub-R-modules of S;
    0 otherwise
minGensIH(I); I=f1,...,fs list of polynomials.
    Returns a minimal  system of generators of <f1,...,fs>, sub-R-module of S
colonInvSyst(f,g); f, g polynomials.
    Returns an element g of R such 	that Aof=g if exists,
    0 otherwise.

MACAULAY INVERSE SYSTEM CORRESPONDENCE WITH NO COEFFICIENTS

invSystGNC(ideal J)
    returns the inverse system of J; J is Artin Gorenstein
invSystNC(J)
    returns the inverse system of J; J is Artin
idealAnnGNC(poly f)
    returns the Artin Gorenstein ideal with inverse system f
idealAnnNC(I)
    returns the Artin  ideal with inverse system I

STRUCTURE OF INJECTIVE HULL WITH NO COEFFICIENTS

memberIHNC(g,I); I=f1,... list of polynomials, g polynomial.
    Returns 1 if g belongs to the R-submodule of S generated by f1,...fs in S,
    returns 0 otherwise
subModIHNC(I,J); I=f1,... list of polynomials, J=g1,... list of polynomials.
    Returns 1 if I is a sub-R-submodule of J, both sub-R-modules of S; 	
    returns 0 otherwise
eqModIHNC(I,J); I=f1,... list of polynomials, J=g1,... list of polynomials.
    Returns 1 if I=J, both sub-R-modules of S;
    returns 0, otherwise
minGensIHNC(I); I=f1,...,fs list of polynomials.
    Returns a minimal system of generators of <f1,...,fs>, sub-R-module of S
colonInvSystNC(f,g); f, g polynomials.
    Returns an element g of R such that Aof=g if exists,
    returns 0 otherwise.

RINGS WITH PARAMETERS

workringc(p, t, n)
    returns the def of  a ring  r with
    t coefficients c	(1),...,c(t),
    n vars x(1),...,x(n), char is p, and order ds

ELLIPTIC CURVES

weiertrassj(t)
    returns the ideal generated by Weierstrass equation
    of the elliptic curve with j invariant
idealwj(t)
    returns the ideal with inverse system weierstrassj(j)
weiertrassp()
    returns the ideal generated by Weierstrass equation of the 	
    elliptic curve with moduli j=c(1)
idealwp(t)
    returns the ideal with inverse system weierstrassp()
    with moduli j=c(1)
\end{verbatim}
}}

\providecommand{\bysame}{\leavevmode\hbox to3em{\hrulefill}\thinspace}
\providecommand{\MR}{\relax\ifhmode\unskip\space\fi MR }
\providecommand{\MRhref}[2]{%
  \href{http://www.ams.org/mathscinet-getitem?mr=#1}{#2}
}
\providecommand{\href}[2]{#2}

\noindent
Juan Elias\\
Departament d'\`Algebra i Geometria\\
Universitat de Barcelona\\
Gran Via 585, 08007 Barcelona, Spain\\
e-mail: {\tt elias@ub.edu}

\end{document}